\documentclass[11pt, a4paper, reqno]{amsart} 
\usepackage[margin=2.8cm]{geometry}
\usepackage{mathabx}
\usepackage[english]{babel}
\usepackage{natbib}
\usepackage{enumitem}
\usepackage[utf8]{inputenc}
\usepackage[dvipsnames]{xcolor}
\usepackage[colorlinks=true, linkcolor=blue, citecolor=blue]{hyperref}
\usepackage[T1]{fontenc}
\usepackage{amsmath}
\usepackage{amsfonts}
\usepackage{amssymb}
\usepackage{amsthm}
\usepackage{bbm}

\patchcmd{\section}{\scshape}{\bfseries}{}{}
\makeatletter
\renewcommand{\@secnumfont}{\bfseries}

\numberwithin{equation}{section}
\newtheorem{thm}{Theorem}[]

\newtheorem{lemma}{Lemma}[section]

\theoremstyle{remark}
\newtheorem*{remark}{Remark}
\newtheorem*{Thm}{Theorem}
\theoremstyle{definition}

\title{Non vanishing of Cubic Twists of $GL_n(\mathbb Q)$ $L$-functions}

\newcommand{\cM}{\mathcal{M}}

\newcommand{\sumstar}{\sideset{}{^*}\sum}
\newcommand{\Lb}{\left(}
\newcommand{\Rb}{\right)}
\newcommand{\md}[1]{\ensuremath{(\operatorname{mod}\, #1)}}

\newcommand\C{\mathbf{C}}

\newcommand{\fp}{\mathfrak{p}}

\newcommand\eps{\varepsilon}

\begin{document}
\author[Sayan Ghosh]{Sayan Ghosh}
	\address{ISI\\Statistics and Mathematics Unit\\
		Kolkata 70010\\
		India}
	\email{gsayan74@gmail.com}
\author[Pratim Mitra]{Pratim Mitra}
	\address{ISI\\Statistics and Mathematics Unit\\
		Kolkata 70010\\
		India}
	\email{pratim2018mitra@gmail.com}
    \date{\today} 

\begin{abstract}
    Let $\pi$ be an irreducible, cuspidal automorphic representation of $GL_n(\mathbb{A}_\mathbb{Q})$ ($n\geq 3$), which is tempered only for $n=3$. Let $s$ be a complex number such that $\Re(s)\notin \left[1/n, 1-1/n\right]$ if $n\neq 4$; $\Re(s)\notin\left[1/5, 4/5\right]$ if $n=4$, then we show that there are infinitely many primitive cubic Dirichlet characters $\chi$ such that $L(s,\pi\times \chi)\neq 0$. Similar results were previously known only for primitive Dirichlet characters without any restriction on the order and quadratic Dirichlet characters.
\end{abstract}
\maketitle{}

\section{Introduction}
One of the central themes of modern analytic number theory is the study of automorphic $L$-functions and their behavior upon twisting by Dirichlet characters. Establishing the non-vanishing theorem for $L$-functions is of fundamental importance, given its profound implications for both automorphic forms and the broader theory of $L$-functions.

In the case of quadratic Dirichlet $L$-functions, Chowla~\cite{Chow65} predicts that central $L$-value, $L(1/2, \chi_d)\neq 0$ for any primitive quadratic Dirichlet character $\chi_d$ of fundamental discriminant $d$. It is now a widely believed folklore that this conjecture holds for all primitive Dirichlet characters. Bui~\cite{Bui12} has shown that at-least $34\%$ of the primitive Dirichlet characters $\chi$ of general moduli $L(1/2, \chi)\neq 0$, and recently this proportion has been improved to $36\%$ by~\cite{QW25}. But for prime modulus, this non-vanishing proportion can be further increased to $38\%$~\cite{KMG22}.  

Non-vanishing for special values of twisted $L$-functions was first studied by Shimura. In the article~\cite{Shim77}, he proved that for a modular form $f$, the twisted $L$-value $L(1/2, f\times\chi)\neq 0$ for infinitely many Dirichlet characters $\chi$. Rohrlich~\cite{Roh89} proved this non-vanishing result in the context of cuspidal automorphic representations of $GL_2(\mathbb A_\mathbb Q)$. This result was first generalized to the higher-rank group $GL(3)$ by Luo~\cite{Luo05}. Recently, building upon his approach, Radziwi\l\l\ and  Yang~\cite{RY24} have further extended this to the limit of analytic methods, ${GL}(4)$. However, to get some non-vanishing results beyond $GL(4)$ using analytic methods, one needs to shift slightly away from the critical line. In this direction, the first contribution was made by Barthel and Ramakrishnan~\cite{BR94}. Although their result is valid for a general number field, we quote it here only for $\mathbb Q$.

Let $\pi$ be an irreducible cuspidal automorphic representation of $GL_n(\mathbb A_\mathbb Q)$ and let $s\in \mathbb C$. Fix any set $S$ of finite places of $\mathbb Q$. Then for $n\geq 3$, there exists infinitely many Dirichlet characters such that $\chi$ is unramified at all places of $S$ and $L(s, \pi\times\chi)\neq 0$ if $\Re(s)\notin \left[\frac{1}{n}, 1-\frac{1}{n}\right]$. This result was further strengthened by Luo~\cite{Luo05}. For $k=3$ the exceptional set is empty and for $n\geq4$, the exceptional set can be replaced by $\left[\frac{2}{n}, 1-\frac{2}{n}\right]$.

A natural extension of this inquiry is to investigate analogous non-vanishing questions for a fixed subfamily of primitive Dirichlet characters. For an integer $l\geq 2$, we define the family of characters of exact order $l$ as,
$$
    \mathcal F_l=\{\chi : \chi\,\, \text{is primitive and of exact order}\,\, l\}.
$$
For $s\in\C$, one can then study the non-vanishing question for the corresponding family of twisted $L$-values,
$$
\{L(s, \pi\times\chi) : \chi\in\mathcal F_l\}.
$$
For quadratic character ($l=2$), this problem has been studied extensively for lower ranks, i.e., for $n=1$ and $2$. When $n=1$, Soundararajan~\cite{Soun00} proved that quadratic Dirichlet $L$-functions $L\left(1/2, \left(\frac{8d}{.}\right)\right)\neq 0$ for at-least $87.5\%$ of fundamental discriminants $8d$ with $d$ odd and $0< 8d\leq X$. For $n=2$  when $\pi$ is attached to a level one modular form, then the recent breakthrough of Li~\cite{Li24} in particularly shows the non-vanishing at the central point. Furthermore for higher ranks $(n\geq 3)$, it was proved by~\cite{CFH05} that there exist infinitely many quadratic characters such that twisted $L$-values as non-zero if $\Re(s)\notin\left[\frac{1}{n}, 1-\frac{1}{n}\right]$. \\ But the non-vanishing problem for $L$-functions twisted by higher-order characters, $l\geq 3$, comparatively less explored. In the case $k=1$, Baier and Young~\cite{BY10} proved the existence of infinitely many cubic Dirichlet characters ($l=3$) such that the $L$-function doesn't vanish at $s=1/2$. Similarly, assuming Lindel\"of hypothesis, Gao and Zhao~\cite{GZ21} proved an analogous result for quartic Dirichlet characters ($l=4$).

The aim of this current article is to investigate this non-vanishing problem in higher rank cases, $n\geq 3$, for the family of primitive cubic Dirichlet characters, which is a thinner and arithmetically subtler family than quadratic characters.

\begin{Thm}
    Let $n\geq 3$, and $\pi$ be an irreducible, cuspidal, automorphic representation of $GL_n(\mathbb A_\mathbb Q)$. Suppose $\pi$ is tempered only for $n=3$, i.e., it satisfies the generalized Ramanujan-Petersson conjecture, then there are infinitely many cubic characters $\chi\in\mathcal F_3$ such that 
    $$L(s, \pi\times\chi)\neq 0,
    \hspace{0.3cm}\text{if}\hspace{0.3cm}\Re(s)\notin
        \begin{cases}
            \hspace{0.1cm}\left[\frac{1}{n}, 1-\frac{1}{n}\right], &{n\neq 4,}\\
           \hspace{0.1cm}\left[\frac{1}{5}, 1-\frac{1}{5}\right], &{n=4}.
        \end{cases}$$
\end{Thm}

Our approach builds on Luo’s in~\cite{Luo05} and uses tools from Baier and Young \cite{BY10}, mainly the large sieve inequality (Lemma \ref{Large_Sieve}). The goal is to obtain an asymptotic identity of the first moment of twisted $L$-values. For $n=3$ instead of averaging over general moduli, we introduce suitable factorization in the moduli. Even though the large sieve inequality used here is imperfect, suitable factorization enables us to use the large sieve \textit{twice} when estimating error term contribution, and it saves more than expected. Other highlighting features are that we haven't used any distribution of the sign of cubic Gau\ss\ sums, and compared to~\cite{BR94} and~\cite{Luo05}, our method is notably elementary in the sense that we haven't used any heavy result, like bounds for hyper-Kloosterman sums. 

Other authors~(\cite{BFH89}, \cite{BFH04}, \cite{CFH05}) have studied the non-vanishing problem for higher order character twists, using the multiple Dirichlet series method. But this method requires the ground field to contain the $l$-th roots of unity. For example, fix $l=3$, and  $\Re(s)\notin\left[\frac{1}{n+1}, \frac{n}{n+1}\right]$. Chinta, Friedberg, and Hoffstein showed that, if the ground field contains \textit{all cube-roots of unity}, then there are infinitely many cubic Hecke characters such that the twisted $L$-value is non-zero. As their method requires, the ground field to contain all cube roots of unity, their method wouldn't work in our case, i.e., for the ground field $\mathbb Q$.
\subsection*{Notations} We are following the $\eps$-convention here, $\eps$ is an arbitrary positive real number, and it may differ at each occurrence. Let $Y>0$ then $X\ll_{\alpha} Y$ or $X=O_\alpha(Y)$ stands for $|X|\leq C Y$ for some absolute constant $C$ depending on some set of parameters $\alpha$. We write $X\asymp Y$ of $Y\ll |X|\ll Y$.

\section{Main Result}
Let $\pi$ be an irreducible, cuspidal representation of $GL_n(\mathbb A_\mathbb Q)$. We will be studying the first moment of the partial $L$- values.
\begin{align}\label{first_moment}
    \mathcal{M}^S_{r_1,r_2}= \sum_{\substack{q_1\; \text{prime}\\ q_1\equiv 1\bmod 3\\ Q^{r_1}\leq q_1\leq 2Q^{r_1}}}\sum_{(q_2,3q_1)=1}W\Lb\frac{q_2}{Q^{r_2}}\Rb\sideset{}{^*}\sum_{\substack{\chi_1\bmod q_1\\ \chi_2\bmod q_2\\\chi_1^3=\chi_2^3=\chi_0}}L^S(s,\pi\times \chi_1\chi_2),
\end{align}
where $^*$ indicates summation taken over primitive characters, $0\leq r_1,r_2\leq1$ with $r_1+r_2=1$, $S$ is a prescribed finite set of primes and $Q$ being large parameter with $Q\rightarrow\infty$. $W$ is a non-negative bump function supported on $[1/2,5/2]$, with $W\equiv 1$ on $[1,2]$.
\begin{thm}\label{Thm1}
    Let $n=3$, and suppose $\pi$ is tempered. Then for every $s$, with $\beta=\Re(s)>2/3$ and $0<r_2<r_1<1$ with $r_1<11/20$ and $r_1+r_2=1$, there exists a set of primes $S$, and $\delta>0$ depending on $\beta,\;r_1,\; r_2$ such that
    \begin{align}\label{asymptotic}
        \mathcal{M}^S_{r_1,r_2}= \cM_{0,0}+O_\delta(Q^{1-\delta}),
    \end{align}
    where $\cM_{0,0}\gg \frac{Q}{\log Q}$. In particular, for $s\in\mathbb C$ such that $\Re(s)>2/3$, there exists infinitely many primitive cubic characters such that, $L(s, \pi\times \chi)\neq0$. 
\end{thm}
\begin{thm}\label{Thm2}
    Let $n\geq 4$, $r_1=0$, and $r_2=1$. Then for every $s$ with $\beta=\Re(s)>\max\{4/5,1-1/n\}$, there exists a set of primes $S$ depending on $n$, and $\delta>0$ depending on $n,\beta$, such that 
        \begin{align}
            \cM^S_{0,1}= \cM_{0,0} + O_\delta(Q^{1-\delta}),
        \end{align}
where $\cM_{0,0}\gg Q$. In particular, for $s\in\mathbb C$ such that $\Re(s)>\max\{4/5,1-1/n\}$, there exists infinitely many primitive cubic characters such that, $L(s, \pi\times \chi)\neq0$.

\begin{remark}
    
\end{remark}
    \begin{itemize}
        \item Our result can be considered as an extension over~\cite{CFH05} except the case $n=4$. 
    \end{itemize}
    \begin{itemize}
        \item For $n=3$, we have chosen a special kind of factorizable moduli, to compensate the ``off-diagonal" contribution~\S\ref{sec: error of M_{1, 2}}. In this process, the factorizable moduli facilitate us to use the cubic large sieve estimate twice instead of once, which saves more than what is expected, and the temperedness of $\pi$ plays a crucial role in this course. 
    \end{itemize}

\end{thm}
\section{Preparatory Tools} 
In this section, we will list a couple of lemmas, which will be used throughout the upcoming sections and are the main tools in our proof.
\subsection{Cubic characters and large Sieve}
Following \cite{BY10} the cubic characters modulo $q$ with $(q,3)=1$, are parametrized by the cubic residue symbol $\chi_\mathfrak n(m)= \Lb\frac{m}{\mathfrak n}\Rb_3$, such that $N(\mathfrak n)=q$  and $\mathfrak n\equiv 1\bmod 3 \in \mathbb{Z}[\omega]$ is square free and has no rational prime divisor. We will use this parametrization for the characters modulo $q_2$ in \eqref{first_moment} to obtain
\begin{align}
    \mathcal M^S_{r_1,r_2}= \sum_{\substack{q_1\; \text{prime}\\ q_1\equiv 1 \bmod 3\\ Q^{r_1}\leq q_1\leq 2Q^{r_1}}}\sideset{}{^*}\sum_{\substack{\chi_1 \bmod q_1\\ \chi_1^3=\chi_0}}\sideset{}{'}\sum_{\substack{\mathfrak n\equiv 1\bmod 3\\ (\mathbf{n},q_1)=1}} L^S(s,\pi \times\chi_1 \chi_{\mathfrak{n}})W\Lb\frac{N(n)}{Q^{r_2}}\Rb,
\end{align}
where the $\sideset{}{'}\sum$ denotes summation taken over $\mathfrak{n}\in \mathbb Z[\omega]$ with $\mathfrak n$ square free and having no rational prime divisor.

Next, we quote the following Large Sieve type inequality, which follows from the proof of Theorem $\text{1.4}$ of \cite{BY10} and the duality principle
\begin{lemma}\label{Large_Sieve}
     Let $a(m)$ be a sequence of complex numbers, supported on an interval of length $\asymp M$. One has,
     \begin{align}
         \sideset{}{'}\sum_{\substack{N(\mathfrak n)\asymp Q\\ \mathfrak n\equiv 1\bmod 3}}\bigg|\sum_{m}a(m)\chi_\mathfrak n(m)\bigg|^2\ll_{\eps} (QM)^\eps(M+Q^{5/3})\Lb\sum_{m}|a(m)|^2\Rb.
     \end{align}
\end{lemma}
\subsection{Second moment for cubic Hecke \texorpdfstring{$L$}{L}-functions}
We will also need the following bound on the second moment of central values of cubic Hecke $L$-functions from $\cite{BY10}$. For $m\in \mathbb{Z}[\omega]$, $\psi_m(\mathfrak n)= \Lb\frac{m}{\mathfrak n}\Rb_3$ defined on the ideals $(\mathfrak n)\subset \mathbb Z[\omega]$ with $\mathfrak n\equiv 1 \bmod 3$, is a Hecke character of modulus $9m$ when $m$ is not an unit or a cube and is the principal character otherwise. The corresponding Hecke $L$-function is defined as
\begin{align}
    L(s,\psi_m)= \sum_{\substack{\mathfrak n\in \mathbb{Z}[\omega]\\ \mathfrak n\equiv 1 \bmod 3}}\Lb\frac{m}{\mathfrak n}\Rb_3N(\mathfrak n)^{-s},
\end{align}
for $\Re(s)>1$. Using the estimate given right after equation (39) of~\cite{BY10}, we have
\begin{lemma}\label{second_moment}
    \begin{align}
        \sum_{m\leq M}\bigg|L\Lb\frac{1}{2}+it,\psi_m\Rb\bigg|^2\ll_\eps M^{3/2+\eps}(1+|t|)^{2/3+\eps}.
    \end{align}
\end{lemma}
\section{Proof of Theorem~\ref{Thm1}}
We will use the approximate functional equation from \cite{Luo05} to express our $L$-values as sums of Dirichlet series.
\subsection{Approximate functional equation}\label{sec: afe} For $f\in C_c^\infty(0,\infty)$ with $\int_0^\infty f(x)\frac{dx}{x}=1,$ set
$$k(w)=\int_0^\infty f(y)y^{w}\frac{dy}{y}.$$
Then $k(w)$ is entire, rapidly decreasing in vertical strips, and $k(0)=1.$ For $x>0$, define
\begin{align}\label{eq: F_1}
    F_1(x)=\frac{1}{2\pi i}\int_{(2)} k(w)x^{-w}\frac{dw}{w},
\end{align}
\begin{align}\label{eq: F_2}
    F_2(x)=\frac{1}{2\pi i}\int_{(2)}k(-w)G(-w+\beta)x^{-w}\frac{dw}{w},
\end{align}
where 
$$G(w)=G(w,\pi)=\frac{L(1-w,\tilde{\pi}_{\infty})}{L(w,\pi_\infty)}.$$
Denote $\beta_0(\pi)=\max_{1\leq j\leq n}\Re(\mu_\infty(j,\pi)).$ Then the functions $F_i(y),$ $i=1,2$ satisy the following:
\begin{enumerate}[label=(\roman*)]
    \item $F_i(y)\ll_n y^{-n}$ for all $n\geq 1$,
    \item $F_1(y)=1+ O(y^n)$ for all $n\geq 1$,
    \item $F_2(y)\ll_{\varepsilon} 1+ y^{1-\beta_0(\pi)-\Re(\beta)-\varepsilon}$,
\end{enumerate}
Applying Cauchy's residue theorem to the inverse Mellin transform ($\sigma>0,X>0$),
\begin{align}
    \frac{1}{2\pi i}\int\limits_{(\sigma)}k(w)L^S(w+s,\pi \times \chi)\Lb\frac{Y}{q^3}\Rb^{-w}\frac{dw}{w}.
\end{align}
one has, for any $X,Y>0$ with $XY= q^k$,
\begin{align}\label{AFE}
    L^S(s,\pi\times \chi)=\sum_{\substack{m=1\\ (m,S)=1}}^\infty \frac{a_\pi(m)\chi(m)}{m^s}F_1\left(\frac{m}{Y}\right)+&\epsilon(s,\pi\times \chi)\sum_{m_1\in \mathcal{I}_S} \frac{b(m_1)}{m_1^s}\\&\sum_{m=1}^\infty\frac{a_{\tilde \pi(m)}\bar \chi(m)}{m^{1-s}}F_2\Lb\frac{m}{Xm_1}\Rb,
\end{align}
where,
\begin{align*}
    \epsilon(s,\pi\times \chi)= W(\pi)\tau(\chi)^3 q^{-3s},
\end{align*}
$\tau(\chi)$ being the Gauss sum
\begin{align*}
    \tau(\chi)=\sum_{r\bmod q}\chi(r)e\Lb\frac{r}{q}\Rb,
\end{align*}
so that $|\tau(\chi)|=q^{1/2}$ and $W(\pi)$ being the root number of $\pi$. The finite Dirichlet polynomial, 
\begin{align}
    \sum_{m_1\in \mathcal{I}_S} \frac{b(m_1)}{m_1^s}= L_S(s,\pi \times \chi)^{-1}= \prod_{p\in S}L_p(s,\pi \times \chi)^{-1},
\end{align}
is the product of the reciprocal of the local factors at $p\in S$. From the classification of the Local representations of $\pi$ at a prime $p$, we see that $\mathcal{I}_S$ contains the set of integers $m_1$ with prime factors only from the set $S$ and such that the largest power of a prime $p$ dividing $m_1$ is at most $k$.
From \eqref{AFE}, we obtain that 
\begin{align}
    \mathcal{M}^S_{r_1,r_2}=\mathcal{M}_1+\mathcal{M}_2,
\end{align}
where,
\begin{align}
    \mathcal{M}_1=&\sum_{\substack{q_1\; \text{prime}\\ q_1\equiv 1 \bmod 3\\ Q^{r_1}\leq q_1\leq 2Q^{r_1}}}\sideset{}{^*}\sum_{\substack{\chi_1\bmod q_1\\ \chi_1^3=\chi_0}} \sideset{}{'}\sum_{\substack {\mathfrak n\equiv 1\bmod 3\\(N(\mathfrak n),q_1)=1}}\\ \nonumber&\sum_{\substack{m=1\\ (m,S)=1}}^\infty \frac{a_{\pi}(m)\chi_1\chi_{\mathfrak{n}}(m)}{m^s} F_1\Lb\frac{mQ^{3/2}}{q^{3}N(\mathfrak n)^3}\Rb W\Lb\frac{N(\mathfrak n)}{Q^{r_2}}\Rb,
\end{align}
and
\begin{align}
    \mathcal{M}_2=&\sum_{\substack{q_1\; \text{prime}\\ q_1\equiv 1 \bmod 3\\ Q^{r_1}\leq q_1\leq 2Q^{r_1}}}\sideset{}{^*}\sum_{\substack{\chi_1\bmod q_1\\ \chi_1^3=\chi_0}}\sideset{}{'}\sum_{\substack{\mathfrak n\equiv 1\bmod 3\\ (N(n),q_1)=1}}\epsilon(s,\pi\times\chi_1 \chi_{\mathfrak{n}})\sum_{m_1\in \mathcal{I}_S} \frac{b(m_1)}{m_1^s}\\ \nonumber&\sum_{m=1}^\infty \frac{a_{\tilde \pi}(m)\widebar{\chi_1\chi_{\mathfrak{n}}}(m)}{m^{1-s}} F_2\Lb\frac{m}{m_1Q^{3/2}}\Rb W\Lb\frac{N(\mathfrak n)}{Q^{r_2}}\Rb.
\end{align}
where $\chi_{\mathfrak{n}}(m)=\Lb\frac{m}{\mathfrak n}\Rb_3$. Before proceeding, we further subdivide $\mathcal{M}_1$ into two sub-components. With $Q^{2r_1/3}\ll Y\ll Q^{5r_1/3}$ we write
\begin{align}
    \mathcal{M}_1= \mathcal{M}_{1,1}+\mathcal{M}_{1,2},
\end{align}
\begin{align}
    \mathcal{M}_{1,1}= &\sum_{\substack{q_1\; \text{prime}\\ q_1\equiv 1 \bmod 3\\ Q^{r_1}\leq q_1\leq 2Q^{r_1}}}\sideset{}{^*} \sum_{\substack{\chi_1\bmod q_1\\ \chi_1^3=\chi_0}}\sideset{}{'}\sum_{\substack {\mathfrak n\equiv 1\bmod 3\\(N(\mathfrak n),q_1)=1}}\\ \nonumber&\sum_{\substack{m=1\\ (m,S)=1}}^\infty \frac{a_{\pi}(m)\chi_1\chi_{\mathfrak{n}}(m)}{m^s} F_1\Lb\frac{mQ^{3r_2}}{YN(\mathfrak n)^3}\Rb W\Lb\frac{N(\mathfrak n)}{Q^{r_2}}\Rb,
\end{align}
and 
\begin{align}
    \mathcal{M}_{1,2}= &\sum_{\substack{q_1\; \text{prime}\\ q_1\equiv 1 \bmod 3\\ Q^{r_1}\leq q_1\leq 2Q^{r_1}}}\sideset{}{^*} \sum_{\substack{\chi_1\bmod q_1\\ \chi_1^3=\chi_0}}\sideset{}{'}\sum_{\substack {\mathfrak n\equiv 1\bmod 3\\(N(\mathfrak n),q_1)=1}}\\ \nonumber&\sum_{\substack{m=1\\ (m,S)=1}}^\infty \frac{a_{\pi}(m)\chi_1\chi_{\mathfrak{n}}(m)}{m^s}\Lb F_1\Lb\frac{mQ^{3r_1}}{Q^{3/2}q_1^3}\Rb- F_1\Lb\frac{mQ^{3r_2}}{YN(\mathfrak n)^3}\Rb\Rb W\Lb\frac{N(\mathfrak n)}{Q^{r_2}}\Rb,
\end{align}
\subsection{The Main Term}\label{sec: main term} We focus our attention to $\mathcal{M}_{1,1}$, which will contribute to the main term. Following \cite{BY10} we express the square freeness and no rational prime divisor condition of $n$ in terms of M\"{o}bius function. Precisely, we get,
\begin{align}
    \mathcal{M}_{1,1}= &\sum_{\substack{q_1\; \text{prime}\\ q_1\equiv 1 \bmod 3\\ Q^{r_1}\leq q_1\leq 2Q^{r_1}}}\sideset{}{^*}\sum_{\substack{\chi\bmod q_1\\ \chi^3=\chi_0}}\sum_{\substack{d\in \mathbb Z\\ d\equiv 1 \bmod 3\\ (d,q_1)=1}}\mu_\mathbb Z(d)\sum_{\substack{l\in \mathbb{Z}[\omega]\\ l\equiv 1\bmod 3\\ (l,q_1)=1}}\mu_\omega(l)\\ \nonumber &\sum_{\substack{m=1\\ (m,S)=1}}^\infty \frac{a_\pi(m)\chi(m)}{m^s}\Lb\frac{m}{dl^2}\Rb_3\mathcal{M}_1(q_1,d,l,m),
\end{align}
where 
\begin{align}
    \mathcal{M}_1(q_1,d,l,m)=\sum_{\substack{n\in \mathbb Z[\omega]\\ n\equiv 1\bmod 3\\ (n,dq_1)=1}}\Lb\frac{m}{n}\Rb_3F_1\Lb\frac{mQ^{3r_2}}{YN(ndl^2)^3}\Rb W\Lb\frac{N(n)}{Q^{r_2}}\Rb,
\end{align}
where $\mu_{\mathbb Z}(d)=\mu(|d|)$, $\mu$ being the usual Mobius function on $\mathbb{N}$ and $\mu_\omega$, the mobius function in $\mathbb{Z}[\omega]$.
As in \cite{BY10}, we use Mellin inversion to express the weight function as a complex integral
\begin{align}
    F_1\Lb\frac{mQ^{r_2n}}{YN(\mathfrak ndl^2)^k}\Rb W\Lb\frac{N(\mathfrak n)}{Q^{r_2}}\Rb= \int\limits_{(2)}^{}\Lb\frac{Q^{r_2}}{N(\mathfrak ndl^2)}\Rb^z\tilde f(z)dz,
\end{align}
where,
\begin{align}
    \tilde f(z)=\int F_1\Lb\frac{m}{Y}x^{-n}\Rb W\Lb x\Rb x^{z-1} dx.
\end{align}
Repeated integration by parts implies that $f$ and its derivatives $\tilde f$ and its derivatives decay like
\begin{align}\label{5.5}
    \ll_j\Lb1 +\frac{m}{Y} \Rb^{-j}\Lb1+|z|\Rb^{-j},\;\; j\geq 0.
\end{align}
So $\mathcal{M}_1(q_1,d,l,m)$ equates to 
\begin{align}
    \mathcal{M}_1(q_1,d,l,m)= \int_{(2)}\left(\frac{Q^{r_2}}{N(dl^2)}\right)^{z}\tilde f(z)L_{dq_1}(z,\psi_m)dz,
\end{align}
where $L_{dq_1}(z,\psi_m)$ denotes the partial Hecke $L$-function of the $\psi_m(n)= \Lb\frac{m}{n}\Rb_3$, with the primes dividing $dq_1$ missing from the Euler product. Shifting the contour to the half line, we encounter a pole at $z=1$ if and only if $m$ is a cube (in which case $\psi_m$ is the principal character). Let $\mathcal{M}_0$ be the contribution of the residues. We then have,
\begin{align}
    \mathcal M_0= Q^{r_2}&\sum_{\substack{q_1\; \text{prime}\\ q_1\equiv 1 \bmod 3\\ Q^{r_1}\leq q_1\leq 2Q^{r_1}}}\sideset{}{^*}\sum_{\substack{\chi_1\bmod q_1\\ \chi_1^3=\chi_0}} \sum_{\substack{d\in \mathbb Z\\ d\equiv 1 \bmod 3\\ (d,q_1)=1}}\frac{\mu_\mathbb Z(d)}{d^2}\sum_{\substack{l\in \mathbb{Z}[\omega]\\ l\equiv 1\bmod 3\\ (l,q_1)=1}}\frac{\mu_\omega(l)}{N(l^2)}\\ \nonumber &\sum_{\substack{m=1\\ (m,dlq_1)=1\\ (m,S)=1}}^\infty\frac{a_\pi(m^3)}{m^{3s}}\tilde{f}(1)\text{Res}_{z=1}L_{dq_1}(z,\psi_m).
\end{align}
Now,
\begin{align*}
    \tilde{f}(1)&= \int F_1\Lb\frac{m^3}{Y}x^{-3}\Rb W(x) dx\\
    &= \frac{1}{2\pi i} \int\limits_{(2)}k(u)\Lb\frac{Y}{m^3}\Rb^u\int W(x)x^{3u}dx \frac{du}{u}\\
    &= \frac{1}{2\pi i}\int\limits_{(2)}k(u)\Lb\frac{Y}{m^3}\Rb^u \widetilde{W}(3u+1)\frac{du}{u}.
\end{align*}
Also,
\begin{align*}
    L_{dq_1}(z,\psi_m)= \zeta_{\mathbb Q(\omega)}(z)\prod_{\fp|3dq_1m}(1-N(\fp)^{-z}), 
\end{align*}
$\zeta_{\mathbb Q(\omega)}(z)$ being the Dedekind zeta function for the field $\mathbb{Q}(\omega)$. Let $c_\omega\ne 0$ be the residue of $\zeta_{\mathbb Q(\omega)}(z)$ at $s=1$. Then
\begin{align}
    \cM_0= 2c_\omega Q^{r_2}&\sum_{\substack{q_1\text{ prime}\\ q_1\equiv 1\bmod 3\\ q_1\sim Q^{r_1}}} \prod_{\fp|q_1}(1-N(\fp)^{-1})\sum_{\substack{m=1\\ (m,q_1)=1\\ (m,S)=1}}^\infty \frac{a_\pi(m^3)}{m^{3s}}\tilde{f}(1)\prod_{\fp|3m}(1-N(\fp)^{-1})\\ \nonumber &\sum_{\substack{d\in \mathbb{Z}\\ d\equiv 1 \bmod 3\\ (d,mq_1)=1}}\frac{\mu_{\mathbb{Z}}(d)}{d^2}\prod_{\fp|d}(1-N(\fp)^{-1})\sum_{\substack{l\in \mathbb{Z}[\omega]\\ l \equiv 1 \bmod 3\\ (l,mq_1)=1}}\frac{\mu_{\omega}(l)}{N(l^2)} .
    \end{align}
Notice that, since $\beta>1/2$,
\begin{align}
    \cM_0= 2c_\omega Q^{r_2}&\sum_{\substack{q_1\text{ prime}\\ q_1\equiv 1\bmod 3\\ q_1\sim Q^{r_1}}} \sum_{\substack{m=1\\(m,S)=1}}^\infty \frac{a_\pi(m^3)}{m^{3s}}\tilde{f}(1)\prod_{\fp|3m}(1-N(\fp)^{-1})\\ \nonumber &\sum_{\substack{d\in \mathbb{Z}\\ d\equiv 1 \bmod 3\\ (d,m)=1}}\frac{\mu_{\mathbb{Z}}(d)}{d^2}\prod_{\fp|d}(1-N(\fp)^{-1})\sum_{\substack{l\in \mathbb{Z}[\omega]\\ l \equiv 1 \bmod 3\\ (l,m)=1}}\frac{\mu_{\omega}(l)}{N(l^2)} +O(Q^{r_2}),
\end{align}
So,
\begin{align}
    \cM_0= {c_\omega} Q^{r_2}\pi(Q^{r_1})&\sum_{\substack{m=1\\ (m,S)=1}}^\infty \frac{a_\pi(m^3)}{m^{3s}}\tilde{f}(1)\prod_{\fp|3m}(1-N(\fp)^{-1})\\ \nonumber &\sum_{\substack{d\in \mathbb{Z}\\ d\equiv 1 \bmod 3\\ (d,m)=1}}\frac{\mu_{\mathbb{Z}}(d)}{d^2}\prod_{\fp|d}(1-N(\fp)^{-1})\sum_{\substack{l\in \mathbb{Z}[\omega]\\ l \equiv 1 \bmod 3\\ (l,m)=1}}\frac{\mu_{\omega}(l)}{N(l^2)} +O(Q^{r_2}).
\end{align}
For $\tilde{f}(1)$, we shift the contour of the $u$- integral to $\sigma=1/3-\beta+\eps$ and collect the residue at $s=0$ to obtain that for a constant $c\neq 0$,
\begin{align}\label{7.11}
    \cM_0= \cM_{0,0} +O(Q^{r_2} + QY^{\frac{1}{3}-\beta+\eps}),
\end{align}
where 
\begin{align*}
    \cM_{0,0} = &Q^{r_2}\Lb\pi(2Q^{r_1},3,1)-\pi(Q^{r_1},3,1)\Rb\sum_{\substack{m=1\\(m,S)=1}}^\infty \frac{a_\pi(m^3)}{m^{3s}}\prod_{\fp|3m}(1-N(\fp)^{-1})\\&\sum_{\substack{d\in \mathbb{Z}\\ d\equiv 1 \bmod 3\\ (d,m)=1}}\frac{\mu_{\mathbb{Z}}(d)}{d^2}\prod_{\fp|d}(1-N(\fp)^{-1})\nonumber \sum_{\substack{l\in \mathbb{Z}[\omega]\\ l \equiv 1 \bmod 3\\ (l,m)=1}}\frac{\mu_{\omega}(l)}{N(l^2)},
\end{align*}
for a constant $c\neq 0$.
Note that if $d<0, d\equiv 1\bmod 3$ then $|d|\equiv 2 \bmod 3$ therefore the $d$-sum reduces to all natural numbers $d$ with $(d,3m)=1$. Since the $l$ sum runs over $\mathbb{Z}[\omega]$ with $l\equiv 1\bmod 3$, $(l,m)=1$, it runs over all non zero the principle ideal in $\mathbb{Z}[\omega]$ which are coprime to $3m$, as each non-zero principal ideal is generated by such an element. Therefore,
\begin{align}
    \cM_{0,0}= &cQ^{r_2}\Lb\pi(2Q^{r_1},3,1)-\pi(Q^{r_1},3,1)\Rb\sum_{\substack{m=1\\(m,S)=1}}^\infty \frac{a_\pi(m^3)}{m^{3s}}\prod_{\fp|3m}(1-N(\fp)^{-1})\\ \nonumber&\hspace{1.25cm}\times \prod_{(p,3m)=1} \Lb1-p^{-2}\Rb\prod_{\fp|p}(1-N(\fp)^{-1})\prod_{(\fp,3m)=1}(1-N(\fp)^{-2}).
\end{align}
Set
\begin{align*}
    \xi= \prod_{p}(1-p^{-2})\prod_{\fp|p}(1-N(\fp)^{-1})\prod_{\fp}(1-N(p)^{-2}).
\end{align*}
Then,
\begin{align}
     \cM_{0,0}&= c\xi Q^{r_2}\Lb\pi(2Q^{r_1},3,1)-\pi(Q^{r_1},3,1)\Rb\sum_{\substack{m=1\\(m,S)=1}}^\infty \frac{a_\pi(m^3)}{m^{3s}}\\ \nonumber&\times\prod_{\fp|3m}(1+N(\fp)^{-1})^{-1}\prod_{p|3m} \Lb1-p^{-2}\Rb^{-1}\prod_{\fp|p}(1-N(\fp)^{-1})^{-1}.
\end{align}
If we specify that $3\in S$, then,
\begin{align}\label{eq:4.26}
    \cM_{0,0}= c'\xi Q^{r_2}&\Lb\pi(2Q^{r_1},3,1)-\pi(Q^{r_1},3,1)\Rb\prod_{p\notin S} \Lb\sum_{j=0}^\infty \frac{a_\pi(p^{3j})}{p^{3js}}\Rb\\& \nonumber\hspace{.7cm}\times\prod_{p\notin S} (1-p^{-2})^{-1}\prod_{\fp | p}(1-N(\fp)^{-2})^{-1},
\end{align}
$c'= \frac{8}{9}c\prod_{\fp | 3}(1-N(\fp)^{-2})^{-1}$.
From the temperedness of $\pi$ one can show that if $S$ contains sufficiently many primes depending on $k$ ,then 
\begin{align*}
    \prod_{p\notin S} \Lb\sum_{j=0}^\infty \frac{a_\pi(p^{3j})}{p^{3js}}\Rb\gg_k 1.
\end{align*}
Hence, for such a choice of the set $S$, the prime number theorem on arithmetic progression implies 
\begin{align*}
    \cM_{0,0}\gg \frac{Q}{\log Q}
\end{align*}

\subsection{Error Contribution from the Shifted Integral}\label{sec: Shifted Integral}

Let $\mathcal{M}_{1,1}'$ be the contribution of the integral shifted to the half-line in the previous section. So,
\begin{align}\label{eq: 4.27}
    \mathcal{M}_{1,1}'&= \sum_{\substack{q_1\; \text{prime}\\ q_1\equiv 1 \bmod 3\\ Q^{r_1}\leq q_1\leq 2Q^{r_1}}}\sideset{}{^*}\sum_{\substack{\chi\bmod q_1\\ \chi^3=\chi_0}}\sum_{\substack{d\in \mathbb Z\\ d\equiv 1 \bmod 3\\ (d,q_1)=1}}\mu_\mathbb Z(d)\sum_{\substack{l\in \mathbb{Z}[\omega]\\ l\equiv 1\bmod 3\\ (l,q_1)=1}}\mu_\omega(l)\\ \nonumber &\times\sum_{\substack{m=1\\ (m,S)=1}}^\infty \frac{a_\pi(m)\chi(m)}{m^s}\Lb\frac{m}{dl^2}\Rb_3\int\limits_{\Lb\frac{1}{2}\Rb}^{}\Lb\frac{Q^{r_2}}{N(dl^2)}\Rb^zL_{dq_1}(z,\psi_m)\tilde f(z)dz.
\end{align}
Using \eqref{5.5}, at the cost of negligible error terms, we can truncate the $m$-sum to $m\ll Y^{1+\eps}$. Similarly, we can truncate the $d$ and $l$-sums to $d\ll Q^{r_2/2+\eps}$, $N(l)\ll Q^{r_2/2+\eps}$. 
One has,
\begin{align*}
    L_{dq_1}(z,\psi_m)= L_d(z,\psi_m)\prod_{\fp |q_1}\Lb1-\frac{\psi_m(\fp)}{N(\fp)^z}\Rb.
\end{align*}
Since $q_1\equiv 1\bmod 3$, upto conjugation, there exists an unique prime $\fp$ in $\mathbb Z[\omega]$ with $N(\fp)=q_1$ such that 
\[q_1=\fp \bar{\fp}.\]
Therefore, fixing such a prime $\fp_{q_1}$ for each $q_1$ and subdividing the $m$-sum into smooth dyadic segments we get,
\[\cM_{1,1}'\ll Q^{\eps}\sup_{1\ll Y_1\ll Y}|E_{0,0}(Y_1)|+|E_{1,0}(Y_1)|+|E_{0,1}(Y_1)|+ |E_{1,1}(Y_1)|\] where
\begin{align}
    E_{i,j}(Y_1)\nonumber&=\sum_{\substack{q_1\; \text{prime}\\ q_1\equiv 1 \bmod 3\\ Q^{r_1}\leq q_1\leq 2Q^{r_1}}}\sideset{}{^*}\sum_{\substack{\chi\bmod q_1\\ \chi^3=\chi_0}}\sum_{\substack{d\in \mathbb Z\\ d\equiv 1 \bmod 3\\ (d,q_1)=1\\ d\ll Q^{r_2/2+\eps}}}\mu_\mathbb Z(d)\\ &\times\sum_{\substack{l\in \mathbb{Z}[\omega]\\ l\equiv 1\bmod 3\\ (l,q_1)=1\\ N(l)\ll Q^{r_2/2+\eps}}}\mu_\omega(l) \sum_{\substack{m\sim Y_1\\ (m,S)=1}}\frac{a_\pi(m)\chi(m)}{m^s}\Lb\frac{m}{dl^2}\Rb_3\\ \nonumber &\times\int\limits_{\Lb\frac{1}{2}\Rb}^{}\Lb\frac{Q^{r_2}}{N(dl^2)}\Rb^z \Lb\frac{\psi_m(\fp_{q_1})}{q_1^z}\Rb^i\Lb\frac{\widebar{\psi_m(\fp_{q_1})}}{q_1^z}\Rb^j L_d(z,\psi_m)\tilde f(z)dz
\end{align}
We first estimate $E_{0,0}(Y_1)$. Rearranging sums and by using the triangle inequality we have, 
\begin{equation*}
    \begin{split}
    E_{0, 0}\ll \underset{d\equiv 1\md{3}}{\sum_{d\ll Q^{r_2/2+\eps}}}\underset{l\equiv 1\md{3}}{\sum_{N(l)\ll Q^{r_2/2+\eps}}}\,\,\underset{(q_1, 3)=1}{\sum_{Q^{r_1}\leq q_1\leq 2Q^{r_1}}}\,\,&\underset{\chi^3=\chi_0}{\sumstar_{\chi\md{q_1}}}\,\,\Bigg|\underset{(m, S)=1}{\sum_{m\sim Y}}\frac{a_\pi(m)\chi(m)}{m^s}\left(\frac{m}{dl^2}\right)_3\\&\times\int\limits_{\Lb\frac{1}{2}\Rb}^{}\Lb\frac{Q^{r_2}}{N(dl^2)}\Rb^z L_d(z,\psi_m)\tilde f(z)dz\Bigg|
    \end{split}
\end{equation*}
Fist we estimate it trivially just by applying Cauchy's inequality on $m\sim Y_1$, and followed by Ramanujan on average bound for Fourier coefficients on average, and the second moment estimate, Lemma \ref{second_moment}. This gives,
\begin{align}\label{E_{0, 0} trivial estiamte}
    E_{00}(Y_1)\ll Q^{\frac{1+r_1}{2}+\eps}Y_1^{\frac{5}{4}-\beta}.
\end{align}
On the other hand, using Cauchy's inequality on the $q_1$-sum, followed by the Large Sieve inequality, i.e., Lemma \ref{Large_Sieve}, we have
\begin{align}
    E_{0,0}(Y_1)\ll Q^{\frac{1}{2}+\eps}  L(Y_1)(Y_1+Q^{\frac{5r_1}{3}})^{\frac{1}{2}},
\end{align}
where 
\begin{align*}
    L(Y_1)= \int \Lb\sum_{m\sim Y_1}\bigg| \frac{a_\pi(m)}{m^s}L_d\Lb\frac{1}{2}+it,\psi_m \Rb\tilde{f}\Lb\frac{1}{2}+it\Rb\bigg|^2\Rb^{1/2}dt,
\end{align*}
where we have estimated the $d$ and $l$-sums trivally in the above estimates.
Again from \eqref{5.5}, Lemma \ref{second_moment} and the temperedness of $\pi$, we obtain
\begin{align}
    L(Y_1)\ll Y_1^{\frac{3}{4}-\beta},
\end{align}
Notice that 
 $Y\ll Q^{5r_1/3}$ forces $Q^{\frac{1}{2}+\eps}Y_1^{\frac{5}{4}-\beta}\ll Q^{\frac{1}{2}+\frac{5r_1}{6}+\eps}Y_1^{\frac{3}{4}-\beta}$ and thus,
\begin{align}
    E_{0,0}\ll Q^{\frac{1}{2}+\eps}\sup_{1\ll Y_1\ll Y}\min\bigg\{Q^{\frac{r_1}{2}}Y_1^{\frac{5}{4}-\beta} , Q^{\frac{5r_1}{6}}Y_1^{\frac{3}{4}-\beta}\bigg\}
\end{align}
For $E_{1,0}$ we obtain, 
\begin{align*}
    E_{1,0}(Y_1) =& \sum_{\substack{q_1\; \text{prime}\\ q_1\equiv 1 \bmod 3\\ Q^{r_1}\leq q_1\leq 2Q^{r_1}}}\sideset{}{^*}\sum_{\substack{\chi\bmod q_1\\ \chi^3=\chi_0}}\sum_{\substack{d\in \mathbb Z\\ d\equiv 1 \bmod 3\\ (d,q_1)=1}}\mu_\mathbb Z(d)\sum_{\substack{l\in \mathbb{Z}[\omega]\\ l\equiv 1\bmod 3\\ (l,q_1)=1}}\mu_\omega(l) \\ \nonumber &\sum_{\substack{m\sim Y_1\\ (m,S)=1}}\frac{a_\pi(m)\chi(m)}{m^s}\Lb\frac{m}{dl^2}\Rb_3\int\limits_{\Lb\frac{1}{2}\Rb}^{}\Lb\frac{Q^{r_2}}{N(dl^2)}\Rb^z\frac{\psi_m(\fp_{q_1})}{q_1^z} L_d(z,\psi_m)\tilde f(z).
\end{align*}
Next, estimating $E_{1,0}(Y_1)$ using Cauchy's inequality, \eqref{5.5}, Lemma \ref{second_moment} and the temperedness of $\pi$ gives,
\begin{align*}
    E_{1,0}(Y_1)\ll Q^{\frac{1}{2}+\eps}Y_1^{\frac{5}{4}-\beta},
\end{align*}
and the same for $E_{0,1},\; E_{1,1}$. Putting everything together,
\begin{align*}
    \cM_{1,1}'\ll_\eps Q^{\frac{1}{2}+\eps}\sup_{1\ll Y_1\ll Y^{1+\eps}} Q^{\frac{r_1}{2}}Y_1^{\frac{3}{4}-\beta}\min\bigg\{Y_1^\frac{1}{2}, Q^{\frac{r_1}{3}}\bigg\}.
\end{align*}
If $\beta \geq 3/4$, the supremum occurs when both the terms inside the curly braces coincide, which occurs at $Y_1\asymp Q^{2r_1/3}$, and we end up with,
\begin{align}\label{error_shifted_integral_1}
    \cM_{1,1}'\ll Q^{\frac{1}{2}+\frac{2r_1(2-\beta)}{3}+\eps}.
\end{align}
Otherwise, the supremum occurs at the extreme range $Y_1\asymp Y$ and
\begin{align}\label{error_shifted_interal_2}
    \cM_{1,1}'\ll_\eps Q^{\frac{1}{2}+\frac{5r_1}{6}+\eps}Y^{\frac{3}{4}-\beta}= Q^{\frac{1}{2}+\eps}\min\bigg\{Q^{\frac{r_1}{2}}Y^{\frac{5}{4}-\beta}, Q^{\frac{5r_1}{6}}Y^{\frac{3}{4}-\beta} \bigg\},
\end{align}
if $Y\gg Q^{2r_1/3}$, which is already satisfied from our choice of Y.
\subsection{Error Contribution from \texorpdfstring{$\cM_{1,2}$}{M{1,2}} and the \texorpdfstring{$\cM_{2}$}{M{2}}}\label{sec: error of M_{1, 2}}
For this section, we assume $r_1>r_2$ (otherwise, can we ensure primitivity of the product $\chi_1\chi_2$). Therefore, we readily have $(q_1,N(n))=1$. Using the decay of $F_1(y)$ for large $y$ and the asymptotic behavior of $F_1$ for small $y$, we can restrict $m$ to the range $Y^{1-\eps}\ll m\ll Q^{3/2+\eps}$ at negligible cost. Therefore, we have,
\begin{align}
    \cM_{1,2}&= \sum_{\substack{q_1\text{ prime}\\ q\equiv 1\bmod 3\\ q_1\sim Q^{r_1}}} \sideset{}{^*}\sum_{\substack{\chi \bmod q_1\\ \chi^3=\chi_0}}\, \sideset{}{'}\sum_{\substack{\mathfrak n\equiv 1\bmod 3}}\,\sum_{\substack{Y^{1-\eps}\ll m\ll Q^{3/2+\eps}\\ (m,S)=1}} \frac{a_\pi(m)\chi_1(m)\chi_{\mathfrak{n}}(m)}{m^s}\\ \nonumber&\hspace{2.1cm}\times\Lb F_1\Lb \frac{mQ^{3r_1}}{Q^{3/2}q_1^{3}}\Rb-F_1\Lb\frac{mQ^{3r_2}}{YN(\mathfrak n)^3}\Rb\Rb W\Lb\frac{N(\mathfrak n)}{Q^{r_2}}\Rb.
\end{align}
Subdividing the $m$-sum into smooth dyadic segments, separating the variables inside the smooth function $F_1$ by its definition (which is a Mellin inversion and results in essentially no cost after shifting the contour to $\Re(w)=\eps$), and taking absolute values inside the $q_1,\chi$-sum, we have,
\begin{align}
   \cM_{1,2}&\ll_\eps  Q^{\eps}\sup_{Y^{1-\eps}\ll Y_1\ll Q^{3/2+\eps}} \sum_{\substack{(q_1,3)=1\\ q_1\sim Q^{r_1}}}\;\sideset{}{^*}\sum_{\substack{\chi \bmod q_1\\ \chi^3=\chi_0}}\\ \nonumber
   &\times\bigg|\sum_{(m,S)=1}\frac{a_\pi(m)}{m^s}G\Lb\frac{m}{Y_1}\Rb\Lb\sideset{}{'}\sum_{\mathfrak n\equiv 1\bmod 3}\chi_{\mathfrak{n}}(m)W\Lb\frac{N(n)}{Q^{r_2}}\Rb\Rb\chi(m)\bigg|,
\end{align}
where $G$ is a compactly supported function on a bounded interval, with derivatives $G^{(j)}\ll_{\eps,j} Q^\eps,\; j\geq 0$.
By Cauchy's inequality, Lemma \ref{Large_Sieve} and the temperedness of $\pi$, 
\begin{align}
    \cM_{1,2}\nonumber&\ll_\eps\sup_{Y^{1-\eps}\ll Y_1\ll Q^{3/2+\eps}} Q^{r_1/2+\eps}Y_1^{-\beta}\Lb Y_1^{\frac{1}{2}}+Q^{\frac{5r_1}{6}}\Rb\\
    &\hspace{.5cm}\times\Lb\sum_{\substack{m\asymp Y_1\\ (m,S)=1}} \bigg|\sideset{}{'}\sum_{\mathfrak n\equiv 1 \bmod 3} \chi_{\mathfrak{n}}(m)W\Lb\frac{N(n)}{Q^{r_2}}\Rb\bigg|^2\Rb^{1/2}.
\end{align}
Finally, along with the duality principle and lemma \ref{Large_Sieve}, it follows that :
\begin{align}
    \cM_{1,2} &\ll_\eps Q^{\frac{1}{2}+\eps}\sup_{Y^{1-\eps}\ll Y_1\ll Q^{3/2}}\Lb Y_1^{1-\beta}+Q^{\frac{5}{6}}Y_1^{-\beta}+ Y_1^{\frac{1}{2}-\beta}Q^{\frac{5r_1}{6}}\Rb,
\end{align}
i.e.
\begin{align}\label{error_M_{1,2}}
    \cM_{1,2}\ll_\eps Q^{\frac{1}{2}+\eps}\Lb Q^{\frac{3(1-\beta)}{2}}+ Q^{\frac{5}{6}}Y^{-\beta}+ Q^{\frac{5r_1}{6}}Y^{\frac{1}{2}-\beta}\Rb.
\end{align}

$\cM_2$ will be estimated similarly. From the decay of $F_2(y)$, the $m$-sum gets restricted to $m\ll m_1Q^{3/2+\eps}$ up to a negligible error. The epsilon factor equals
\[\epsilon(s,\pi\times\chi\chi_{\mathfrak{n}})= W(\pi) \tau(\chi_1)^3\tau(\chi_{\mathfrak{n}})^3q_1^{-3s}N(n)^{-3s}.\] 
Observe that for the finite Dirichlet polynomial, 
\begin{align*}
    L_S(s,\pi\times \chi)^{-1}= \sum_{m_1\in \mathcal{I}_S}\frac{b(m_1)}{m_1^s}, 
\end{align*}
one has 
\begin{equation}\label{eq:b(m_1)}
    b(m_1)= c(m_1)\chi(m_1),
\end{equation} for a complex number $c(m_1)$ depending only on the local representation of $\pi$ at the primes dividing $m_1$. For example, if the local representations are principal series and if $\chi$ stays unramified at the corresponding primes, then the coefficient $c(m_1)$ is a finite linear combination of the powers of local Satake parameters. 
Again, subdividing the $m$-sum into smooth dyadic segments, we have,
\begin{align}
    \cM_2&\ll_{k,\eps} Q^\eps|\mathcal{I}_S| \sup_{\substack{m_1\in \mathcal{I}_S\\ 1\ll Y_1\ll m_1Q^{3/2+\eps}}}|c(m_1)|\sum_{\substack{(q_1,3)=1\\ q_1\asymp Q^{r_1}}}q_1^{3\Lb\frac{1}{2}-\beta\Rb}\sideset{}{^*}\sum_{\substack{\chi \bmod q_1\\ \chi^3=\chi_0}}\\ \nonumber&\times \bigg|\sum_{m} \frac{a_{\tilde{\pi}}(m)}{m^{1-s}}G\Lb\frac{m}{Y_1}\Rb\Lb \sideset{}{'}\sum_{\mathfrak n\equiv 1\bmod 3} \frac{\tau(\chi_{\mathfrak{n}})^3\chi_{\mathfrak{n}}(m_1)}{N(n)^{3s}}\bar{\chi}_n(m) W\Lb\frac{N(n)}{Q^{r_2}}\Rb\Rb\chi(m)\bigg|. 
\end{align} 
Estimating similarly as \eqref{error_M_{1,2}} we obtain
\begin{align}
    \cM_2 \ll_\eps Q^{3\Lb\frac{1}{2}-\beta\Rb+\frac{1}{2}+\eps}\sup_{1\ll Y_1\ll Q^{\frac{3}{2}}} \Lb Y_1^\beta+ Q^{\frac{5}{6}}Y_1^{\beta-1}+ Q^{\frac{5r_1}{6}}Y_1^{\beta-\frac{1}{2}}\Rb,
\end{align}
i.e.
\begin{align}\label{error_M_2}
    \cM_2 \ll_\eps Q^{3\left(\frac{1}{2}-\beta\right)+\frac{1}{2}+\eps}\Lb Q^{\frac{5}{6}}+ Q^{\frac{3\beta}{2}}+ Q^{\frac{5r_1}{6}+\frac{3\left(\beta-1/2\right)}{2}}\Rb
\end{align}
\subsection{Conclusion}
Fixing $r_1,\;r_2$ such that $9/20<r_2<r_1<11/20$ and $r_1+r_2=1$, we compare the terms from \eqref{error_shifted_integral_1}, \eqref{error_shifted_interal_2}, \eqref{error_M_{1,2}} and \eqref{error_M_2} getting,
\begin{equation}\label{4.44}
    \cM_{1,1}'+\cM'_{1,2} +\cM_2 \ll_\eps Q^{\frac{1}{2}+\eps}\Lb Q^{\frac{5r_1}{6}}Y^{\frac{3}{4}-\beta}+ Q^{\frac{5}{6}}Y^{-\beta} +Q^{\frac{3}{2}(1-\beta)}+ Q^{\frac{5}{6}-3\left(\beta-\frac{1}{2}\right)} \Rb,
\end{equation}
if $\beta<3/4$ and
\begin{align}\label{10.2}
    \cM_{1,1}'+\cM'_{1,2} +\cM_2 \ll_\eps Q^{\frac{1}{2}+\eps}&\Lb Q^{\frac{2r_1(2-\beta)}{3}}+ Q^{\frac{3(1-\beta)}{2}}+Q^{\frac{5}{6}}Y^{-\beta}+ Q^{\frac{5r_1}{6}}Y^{\frac{1}{2}-\beta} + Q^{\frac{7}{3}-3\beta}\Rb,
\end{align}
if $\beta\geq 3/4$.\\
Choosing any $Y\in \Lb Q^{1/2},Q^{2(3-5r_1)}\Rb$ if $2/3<\beta<3/4$, and $Y\in \Lb Q^{4/9}, Q^{5r_1/3} \Rb$ if 
$\beta\geq 3/4$, we find by checking each term separately that there is  $\delta=\delta(\beta,r_1,r_2)>0$, such that 
\begin{align*}
    \cM_{1,1}'+\cM'_{1,2} +\cM_2 \ll_\delta Q^{1-\delta}.
\end{align*}
Together with \eqref{7.11} we obtain
\begin{align}
    \cM^S_{r_1r_2}= \cM_{0,0}+ O_\delta(Q^{1-\delta})
\end{align}
for some $\delta=\delta(\beta,r_1,r_2)>0$, which culminates our proof of Theorem~\ref{Thm1}.\qed

\section{Proof of Theorem~\ref{Thm2}}
In this section, we are going to prove the non-vanishing result for the remaining case, $n\geq 4$. Our proof will be almost similar to the proof of case $n=3$ covered by Theorem~\ref{Thm1}. So, to avoid unnecessary computational details, we briefly outline the key steps below.\\

Like before~\S\ref{sec: afe} we begin by using the approximate functional equation to rewrite $\cM_{0, 1}$ as,
$$\cM_{0, 1}^S=\cM_1+\cM_2,$$
where
\begin{align}
    \cM_1=\sideset{}{'}\sum_{\substack{\mathfrak n\equiv 1\bmod 3\\}}\sum_{\substack{m=1\\(m,S)=1}}^\infty \frac{a_{ \pi}(m){\chi_{\mathfrak{n}}}(m)}{m^{s}} F_1\Lb\frac{mQ^{n/2}}{N(\mathfrak n)^n}\Rb W\Lb\frac{N(\mathfrak n)}{Q}\Rb,
\end{align}
and
\begin{align}
    \cM_2=\sideset{}{'}\sum_{\substack{\mathfrak n\equiv 1\bmod 3\\}}\epsilon(s,\pi\times \chi_{\mathfrak{n}})\sum_{m_1\in \mathcal{I}_S} \frac{b(m_1)}{m_1^s}\sum_{m=1}^\infty \frac{a_{\tilde \pi}(m)\widebar{\chi_{\mathfrak{n}}}(m)}{m^{1-s}} F_2\Lb\frac{m}{m_1Q^{n/2}}\Rb W\Lb\frac{N(\mathfrak n)}{Q}\Rb,
\end{align}
where $F_i$'s are given by~\eqref{eq: F_1}, and~\eqref{eq: F_2}. Now for some parameter $1\ll Y\ll Q^{n/2}$ (to be chosen later), we further subdivide $\cM_1$ as \\$\cM_1=\cM_{1,1}+\cM_{1,2}$, where
\begin{align}
    \cM_{1,1}=\sideset{}{'}\sum_{\substack{\mathfrak n\equiv 1\bmod 3\\}}\sum_{\substack{m=1\\(m,S)=1}}^\infty \frac{a_{ \pi}(m){\chi_{\mathfrak{n}}}(m)}{m^{s}} F_1\Lb\frac{m}{Y}\Rb W\Lb\frac{N(\mathfrak n)}{Q}\Rb,
\end{align}
and
\begin{align}
    \cM_{1,2}=\sideset{}{'}\sum_{\substack{\mathfrak n\equiv 1\bmod 3\\}}\sum_{\substack{m=1\\(m,S)=1}}^\infty \frac{a_{ \pi}(m){\chi_{\mathfrak{n}}}(m)}{m^{s}} \Lb F_1\Lb\frac{mQ^{n/2}}{N(\mathfrak n)^n}\Rb - F_1\Lb\frac{m}{Y}\Rb\Rb W\Lb\frac{N(\mathfrak n)}{Q}\Rb.
\end{align}
We extract our main term from $\cM_{1,1}$ by estimating it exactly as \S\ref{sec: main term} and \S\ref{sec: Shifted Integral}. Precisely, writing the weight functions in terms of Mellin inversions, shifting the contours to the left and collecting the residues at $z=1$ if $m$ is a cube, obtaining
\begin{align}\label{second main}
    \cM_{1,1}=\cM_{0}+ \cM'_{1,1},
\end{align}
\begin{align*}
    \mathcal M_0= Q&\sum_{\substack{d\in \mathbb Z\\ d\equiv 1 \bmod 3\\}}\frac{\mu_\mathbb Z(d)}{d^2}\sum_{\substack{l\in \mathbb{Z}[\omega]\\ l\equiv 1\bmod 3\\ }}\frac{\mu_\omega(l)}{N(l^2)}\\ \nonumber &\sum_{\substack{m=1\\ (m,dl)=1\\ (m,S)=1}}^\infty\frac{a_\pi(m^3)}{m^{3s}}\tilde{f}(1)\text{Res}_{z=1}L_{d}(z,\psi_m)
\end{align*}
and $\cM'_{1,1}$ is the contribution of the ``shifted integral" to $\Re(z)=1/2$. Evaluating $\cM_0$ in a similar fashion as \eqref{7.11} we find
\begin{align}\label{5.6}
    \cM_{0}= \cM_{0,0}+O(QY^{\frac{1}{3}-\beta+\eps}),
\end{align}
where 
\begin{align*}
    \cM_{0,0}=cQ \prod_{p\notin S}\Lb\sum_{j=0}^\infty \frac{a_\pi(p^{3j})}{p^{3js}}\Rb,
\end{align*}
$c$ is a non-zero constant depending on $s$, which can be expressed as an Euler product similar to~\eqref{eq:4.26}. By the Jacquet-Shalika bound for the Fourier coefficients, $a_\pi(m)\ll m^{1/2}$, since $\beta>1/2$, we see that 
$\cM_{0,0}\gg_n Q$, if $S$ is contains sufficiently many primes depending on $n$.

We express $\cM'_{1,1}$ and truncate the $d$, $l$, $m$ sums to $d\ll Q^{1/2+\eps}$, $N(l)\ll Q^{1/2+\eps}$ and $m\ll Y^{1+\eps}$ as as \eqref{eq: 4.27}. Next, taking absolute value inside the $d,~l$ sums and estimating the $m$ sum by Cauchy's inequality, Ramanujan bound on average, and the second moment bound in Lemma \ref{second_moment}, we get
\begin{align}\label{5.7}
    \cM_{1,1}'\ll_\eps Q^{\frac{1}{2}+\eps} Y^{\frac{5}{4}-\beta}.
\end{align}
Since we are not introducing any factorization of the moduli here, there are no harmonics left to apply the Large Sieve inequality as we did before.
Plugging the estimates~\eqref{5.6}, and~\eqref{5.7} in~\eqref{second main} we get, 
\begin{equation}\label{main term}
    \cM_{1, 1}= cQ \prod_{p\notin S}\Lb\sum_{j=0}^\infty \frac{a_\pi(p^{3j})}{p^{3js}}\Rb +O_\eps(QY^{\frac{1}{3}-\beta+\eps}+Q^{\frac{1}{2}+\eps} Y^{\frac{5}{4}-\beta}).
\end{equation}
Hence, we are done with the evaluation of the main term. 

Next, we will estimate the error term contribution coming from $\cM_{1,2}$ and $\cM_{2}$. Both of them will be estimated similarly along the lines of \S \ref{sec: error of M_{1, 2}}, so we just outline the proof of $\cM_{2}$. We have $\epsilon(\pi \times \chi_\mathfrak n,s)=W(\pi)\tau(\chi_\mathfrak n)^nN(\mathfrak n)^{-ns}\asymp Q^{n(1/2-\beta)}$. We take absolute value inside the $m_1$ and $\mathfrak{n}$ sums getting 
\begin{align*}
    \cM_2= Q^{n\Lb \frac{1}{2}-\beta \Rb}\sup_{m_1\in \mathcal{I}_S} \sideset{}{'}\sum_{\mathfrak n\equiv 1\md{3}} \sum_{m=1}^\infty\bigg| \frac{a_{\tilde \pi}(m)\widebar{\chi_{\mathfrak{n}}}(m)}{m^{1-s}} F_2\Lb\frac{m}{m_1Q^{n/2}}\Rb\bigg|
\end{align*}
Subdividing the $m$ sums into smooth dyadic segments of size $1\ll Y_1\ll m_1Q^{n/2+\eps}$ 
By using the Cauchy's inequality over $\mathfrak n$ and $m_1$-sum we have, 
\begin{equation*}
    \cM_2\ll_\eps Q^{\frac{n+1}{2}-n\beta+\eps}\sup_{\substack{m_1\in \mathcal{I}_S\\ 1\ll Y_1\ll m_1Q^{n/2+\eps}}}\Lb\sideset{}{'}\sum_{\substack{\mathfrak n\equiv 1\md{3}\\ N(\mathfrak n)\asymp Q}}\bigg|\sum_{m\asymp Y_1} \frac{a_{\tilde \pi}(m)\widebar{\chi_{\mathfrak{n}}}(m)}{m^{1-s}} F_2\Lb\frac{m}{m_1Q^{n/2}}\Rb\bigg|^2\Rb^{1/2}
\end{equation*}

Now using the cubic large sieve estimate~(Lemma~\ref{Large_Sieve}) we have, 

\begin{align}\label{error: M_2}
    \cM_2\ll_\eps Q^{\frac{n+1}{2}-n\beta+\eps}\sup_{\substack{m_1\in\mathcal I_S\\1\ll Y_1\ll Q^{n/2+\eps}}}(Y_1^{\frac{1}{2}}+Q^{\frac{5}{6}})Y_1^{\beta-\frac{1}{2}}\ll_\eps Q^{\frac{n+1}{2}-n\beta+\eps}(Q^{\frac{5}{6}}+ Q^{\frac{n\beta}{2}}).
\end{align}
Along the same lines, we obtain,
\begin{align}\label{error: M_{1, 2}}
    \cM_{1,2}\ll_{\eps} \sup_{Y^{1-\eps}\ll Y_1\ll Q^{n/2}}Q^{\frac{1}{2}+\eps}Y_1^{\frac{1}{2}-\beta}(Y_1^{\frac{1}{2}}+Q^{\frac{5}{6}})\ll_\eps Q^{\frac{1}{2}+\eps}\Lb Q^{\frac{n(1-\beta)}{2}}+ Q^{\frac{5}{6}}Y^{\frac{1}{2}-\beta}\Rb.
\end{align}
If $\beta>\max\{4/5,1-1/n\}$, then one can choose $Y\in \Lb Q^{\frac{2}{3(2\beta-1)}}, Q^{\frac{2}{5-4\beta}}\Rb$ so that \eqref{error: M_2}, \eqref{error: M_{1, 2}} and the error term of \eqref{main term} are all $O_\delta (Q^{1-\delta}),$ for some $\delta=\delta(n,\beta)>0$. This concludes the proof of Theorem~\ref{Thm2}.\qed
\subsection*{Acknowledgment} The authors would like to thank Prof. Ritabrata Munshi for various discussions and valuable suggestions. We are also grateful to the Indian Statistical Institute, Kolkata, for providing an excellent and productive research environment.
\bibliographystyle{plain}
\bibliography{ref}

\end{document}